\documentclass[11pt]{article}
\usepackage{latexsym,amssymb,amsmath,amscd,amsthm,amsxtra}
\usepackage{mathrsfs}
\usepackage{epsfig}
\usepackage{pgf,tikz}
\usepackage{graphicx}
\usepackage{graphics}
\usepackage{setspace}
\newcommand{\counte}{section}

\newtheorem{prop}{\bf Proposition}[\counte]
\newtheorem{lemma}{\bf Lemma}[\counte]

\newtheorem{coro}{\bf Corollary}[\counte]

\newtheorem{remark}{\bf Remark}[\counte]

\newcommand{\D}{\displaystyle}

\author{Zhao Xu-an, zhaoxa@bnu.edu.cn\\ Gao Hongzhu, hzgao@bnu.edu.cn\\Department of Mathematics, Beijing Normal University\\Key Laboratory
of Mathematics and Complex Systems\\ Ministry of Education,
China, Beijing 100875}

\title{The rational cohomology Hopf algebra of a generic Kac-Moody group\thanks{The authors are supported by National Science Foundation of China, 11571038.}}

\date{}

\begin{document}

\maketitle
\begin{abstract}
In this paper we determine the rational homotopy type of the classifying space of a generic Kac-Moody group by computing its rational cohomology ring. As an application we determine the rational homology Hopf algebra of the generic Kac-Moody group.
\end{abstract}



\section{Introduction}

Let $A=(a_{ij})$ be an $n\times n$ Cartan matrix. By Kac\cite{Kac_68} and Moody\cite{Moody_68}, it is well known that there is a Kac-Moody Lie algebra $g(A)$ associated to $A$. In \cite{Kac_Peterson_83}\cite{Kac_Peterson_84}\cite{Kac_85} Kac and Peterson constructed the corresponding Kac-Moody group $G(A)$. In this paper for convenience we consider the derived Lie algebra $g'(A)$ and the associated simply connected group $G'(A)$. But we still use the symbols $g(A)$ and $G(A)$.



Cartan matrices are divided into three types, i.e. finite type, affine type and indefinite type. A Cartan matrix $A$ is indecomposable if $A$ can't be written as a direct sum of two Cartan matrices $A_1$ and $A_2$. $A$ is symmetrizable if there exists an invertible diagonal matrix $D$ and a symmetric matrix $B$ such that $A=DB$, see Kac\cite{Kac_82} for details. A Cartan matrix $A$ is generic if $a_{ij}a_{ji}\geq 4$ for all $i,j$. $A$ is generic if and only if all its principal sub-matrices of rank $2$ are not of finite type. All these properties for Cartan matrices  can be used for the associated Kac-Moody Lie algebras and Kac-Moody groups. For example a generic Kac-Moody group is indecomposable.

In \cite{Zhao_Jin_15}, the authors determined the rational homotopy type of the indefinite Kac-Moody group $G(A)$. Since $G(A)$ has a multiplication, it is a Hopf space. It is important to determine the rational Hopf homotopy type of $G(A)$. This is equivalent to determine the rational cohomology Hopf algebra $H^*(G(A))$ or the dual rational homology Hopf algebra $H_*(G(A))$. It is further equivalent to determine the rational cohomology algebra $H^*(BG(A))$ of the classifying space $BG(A)$.

On the rational homotopy group $\pi_*(G(A))$, the Samelson product $[\  , \ ]:\pi_p(G(A))\times \pi_q(G(A))\to \pi_{p+q}(G(A))$ is defined as
$$[\alpha,\beta](s\wedge t)=\alpha(s)\beta(t)\alpha(s)^{-1}\beta(t)^{-1},s\in S^p,t\in S^q.$$
$(\pi_*(G(A)),[\  , \ ])$ is a rational graded Lie algebra.



Let $\chi: \pi_*(G(A))\to H_*(G(A))$ be the Hurewicz morphism of graded Lie algebras. By Milnor and Moore\cite{Milnor_Moore_65}, the induced morphism $\widetilde \chi: U(\pi_*(G(A)))\to H_*(G(A))$ is an isomorphism of Hopf algebras, where $U(\pi_*(G(A)))$ is the universal enveloping algebra of $\pi_*(G(A))$. And $H_*(G(A))$ is primitively generated by $\pi_*(G(A))$.
So to determine the Hopf algebra structure on $H_*(G(A))$, it is enough to compute the graded Lie algebra $\pi_*(G(A))$. By combining rational homotopy theory(see \cite{Sullivan_77}) with the cohomology ring $H^*(G(A))$(see \cite{Zhao_Jin_13}), one knows that for a generic Cartan matrix $A$, $\pi_{odd}(G(A))\cong \mathbb{Q}$ or $\{0\}$ depending on whether $A$ is symmetrizable or not. And the decomposition $\pi_*(G(A))=\pi_{even}(G(A))\oplus \pi_{odd}(G(A))$ is a decomposition of Lie algebras. $H_{even}(G(A))$(i.e. the rational Chow ring of $G(A)$) is isomorphic to the universal enveloping algebra  $U(\pi_{even}(G(A)))
$. By \cite{Kac_85} the Poincar\'{e} series of $H_{even}(G(A))$ is $$C_A(q)=P_{F(A)}(q)(1-q^2)^n(1-q^4)^{-\epsilon(A)}.$$
Here $P_{F(A)}(q)$ is the Poincar\'{e} series of the flag manifolds $F(A)$ associate to Kac-Moody group $G(A)$ and $\epsilon(A)=\dim \pi_{odd}(G(A))$. By \cite{Zhao_Jin_13}, $\D{P_{F(A)}(q)=\frac{1+q^2}{1-(n-1)q^2}}$. Hence we can compute the Poincar\'{e} series of $H_{even}(G(A))$.

For a non-symmetrizable Cartan matrix $A$, we write the Poincar\'{e} series $\D{C_A(q)=\frac{1}{\frac{1-(n-1)q^2}{(1-q^2)^{n-1} (1-q^4)}}}$ as
$$\D{\frac{1}{1-a_4 q^4-a_4 q^6-\cdots-a_{2i} q^{2i}-\cdots}},$$ and for a symmetrizable Cartan matrix $A$, we write $\D{C_A(q)=\frac{1}{\frac{1-(n-1)q^2}{(1-q^2)^{n-1}}}}$ as $$\D\frac{1}{1-b_4 q^4-b_6 q^6-\cdots-b_{2i} q^{2i}-\cdots},$$
where for $i\geq 2$, $a_{2i}=\sum\limits_{k=0}^{[\frac{i}{2}]}(i-1-2k){{n+i-2k-3} \choose{n-3}}, b_{2i}= (i-1){{n+i-3}\choose{n-3}}$ are natural numbers depending on $n$.

These two Poincar\'{e} series are the same as the Poincar\'{e} series of the tensor Hopf algebra with $a_{2i}$ and $b_{2i}$ generators of degree $2i$ for each $i\geq 2$. In \cite{Zhao_Jin_15} we gave the following conjecture.

\noindent{\bf Conjecture: }For a non-symmetrizable(or symmetrizable) generic Kac-Moody group $G(A)$, the graded Lie algebra $\pi_{even}(G(A))$ is a free Lie algebra with $a_{2i}$(or $b_{2i}$) generators of degree $2i$ for each $i\geq 2$.

Since the universal enveloping algebra of an even graded free Lie algebra is a tensor algebra, if the conjecture is true, then $H_{even}(G(A))$ is a tensor Hopf algebra with $a_{2i}$(or $b_{2i}$) generators of degree $2i$ for each $i\geq 2$.



In this paper we compute the rational cohomology ring $H^*(BG(A))$ at first. This determines the rational homotopy type of $BG(A)$. Then we compute the graded Lie algebra $\pi_{*}(BG(A))$ with Whitehead product. Since the graded Lie algebra $\pi_{*}(G(A))$ with Samelson product is determined by $\pi_*(BG(A))$ with Whitehead product, we determine the graded Lie algebra $\pi_*(G(A))$ and prove the conjecture at last.

The main results in this paper are the following two theorems.

\noindent{\bf Theorem 1:}
If $A=(a_{ij})_{n\times n}$ is a non-symmetrizable generic Cartan matrix, $n\geq 3$, 
then the Poincar\'{e} series of $BG(A)$ is $\D{P_n(q)=q\left [\frac{(n-1)q^2-1}{(1-q^2)^{n-1} (1-q^4)} +1\right ]+1}$.

\noindent{\bf Theorem 2:}
If $A=(a_{ij})_{n\times n}$ is a symmetrizable generic Cartan matrix, $n\geq 2$, then the Poincar\'{e} series of $BG(A)$ is $\D{Q_n(q)=\frac{1}{1-q^4}(q\left [ \frac{(n-1)q^2-1}{(1-q^2)^{n-1}}  +1\right ]+1)}$.

The contents of this paper are as follows: in section 2 we give some preparing lemmas, in section 3 we prove Theorem 1 and 2, in section 4 we give some results derived from these theorems, including the above conjecture.

\section{Some preparing lemmas}
In the following all the Cartan matrices are assumed to be generic. All the homology and cohomology are of rational coefficients.

Let $S$ be the set of integers $1,2,\cdots,n$, and $\Pi=\{\alpha_1,\alpha_2,\cdots,\alpha_n\}$ be the simple root system of $G(A)$. For each $I\subset S$, the matrix $A_I=(a_{ij})_{i,j\in I}$ is also a Cartan matrix. Corresponding to $I\subset S$, there is a parabolic subgroup $G_I(A)$ of $G(A)$ whose simple root system is $\Pi_I=\{\alpha_i|i\in I\}$. All the proper subsets of $S$ form a category $\mathrm{C}$ with object $I\subset S$ and morphism $I\subset J$. By constructing classifying spaces we have a functor $F: \mathrm{C}\to Top$ which sends $I$ to $BG_I(A)$ and $I\subset J$ to map $BG_I(A)\to BG_J(A)$.

Since we only consider the homotopy type of the Kac-Moody group we replace the group $G(A)$(or $G_I(A)$) by its unitary form and use the same symbol.

We need the following lemmas to prove the main theorems.

\begin{lemma}
For a Kac-Moody group $G(A)$ and $I\subset S$, the subgroup $G_I(A)$ is isomorphic to $G(A_I)\widetilde \times T^{n-|I|}$, the semi-direct product of $G(A_I)$ and $T^{n-|I|}$. As a result there is an isomorphism $H^*(BG_I(A))\cong H^*(BG(A_I))\times H^*(BT^{n-|I|})$.
\end{lemma}

By this lemma, the Poincar\'{e} series of $BG_I(A)$ is obtained from the Poincar\'{e} series of $BG(A_I)$ by multiplying a factor $\D{\frac{1}{(1-q^2)^{n-|I|}}}$.

By Kichiloo\cite{Kitchloo_98}\cite{Broto_Kitchloo_02}, for a Cartan matrix of infinite type, the homotopy colimit of the functor $F$ gives the homotopy type of $BG(A)$. For any $I\in \mathrm{C}$, let $\mathrm{C}_I$ be the full subcategory of $\mathrm{C}$ whose objects are proper subsets of $I$. If $|I|\geq 2$, then $G_I(A)$
is of infinite type. By using the result of Kichiloo to $BG_I(A)$, we get $H^*(BG_I(A))\simeq \textrm{colimit} F|_{\mathrm{C}_I}$. As a consequence we have

\begin{lemma}
Let $\mathrm{C}'$ be the full subcategory of $\mathrm{C}$ which contains only objects $\emptyset,\{1\},\{2\},\cdots,\{n\}$, then for a generic Kac-Moody group $G(A)$,

$BG(A)\simeq \mathrm{colimit}F|_{\mathrm{C}'}$

$\ \ \ \ \ \ \ \ \  \  \simeq BG_1(A)\cup_{BT} BG_2(A)\cup_{BT}\cdots \cup_{BT} BG_n(A)$

$\ \ \ \ \ \ \ \ \  \ \simeq BG_{\{1,2,\cdots,n-1\}}(A)\cup_{BT} BG_n(A)$.
\end{lemma}




The action of Weyl group $W(A)$(or $W_I(A)$) of $G(A)$(or $G_I(A)$) on the maximal torus $T$ induces the action of $W(A)$(or $W_I(A)$) on $H^*(BT)$.

\begin{lemma}
For a generic Kac-Moody group $G(A)$, the image of the homomorphism $Bi_I^*: H^*(BG_I(A))\to H^*(BT)$ induced by the inclusion $i_I:T\subset G_I(A)$ is $H^*(BT)^{W_I(A)}$, i.e. the $W_I(A)$ invariants. In particularly the image of the homomorphism $H^*(BG(A))\to H^*(BT)$ is $H^*(BT)^{W(A)}$.
\end{lemma}

This lemma is the generalization of a result of Borel\cite{Borel_53_1} for compact Lie groups. It can be proved in the inductive procedure of the proofs for the main theorems. see \cite{Zhao_Gao_Ruan_18} for details.

\begin{lemma}
If $A$ is a non-symmetrizable generic Cartan matrix, then there exist $i,j,k\in S,i<j<k$ such that $A_{\{i,j,k\}}$ is non-symmetrizable.
\end{lemma}

\noindent {\bf Proof: }Suppose this lemma is not true for Cartan matrix $A$. Then for any $i<j<k$, $A_{\{i,j,k\}}$ is symmetrizable. Hence $a_{ij}a_{jk}a_{ki}=a_{ik}a_{kj}a_{ji}$. We set $d_{i}=\D{\frac{a_{1i}}{a_{i1}}}$. For $1<i<j$, we have $a_{1i}a_{ij}a_{j1}=a_{1j}a_{ji}a_{i1}$, hence $\D{\frac{a_{ij}}{a_{ji}}=\frac{d_j}{d_i}}$. But this means that $A$ is symmetrizable, a contradiction.

\begin{lemma}
Let $A$ be a generic Cartan matrix. If $A$ is symmetrizable, then $H^*(BT)^{W(A)}\cong \mathbb{Q}[\psi]$, where $\psi$ is the Killing form; If $A$ is non-symmetrizable, then $H^*(BT)^{W(A)}\cong\mathbb{Q}$.
\end{lemma}

This result was proved in Zhao-Jin\cite{Zhao_Jin_12}. In fact it is valid for an arbitrary indefinite and indecomposable Cartan matrix.

\begin{lemma}
Let $X=X_1\cup_{X_0} X_2$ be the push-out of the diagram $X_1 \stackrel{j_1} \leftarrow X_0\stackrel{j_2}\rightarrow X_2$. The homomorphism $j: H^*(X_1)\oplus H^*(X_2)\to H^*(X_0)$ is given by $j(u,v)=j_1^*(u)-j_2^*(v)$. If $X_1,X_2$ are deformation retracts of some open subspaces of $X$, then there exists a short exact sequence $$0\to \Sigma \mathrm{coker}j\to H^*(X)\to \ker j\to 0.$$
\end{lemma}

\noindent{\bf Proof: }We have the Mayer-Vietoris exact sequence
$$\cdots\to H^{*-1}(X_{1})\oplus H^{*-1}(X_{2})\stackrel{j}\to H^{*-1}(X_0)\stackrel{\delta}\to H^*(X)\stackrel{i}\to H^*(X_{1})\oplus H^*(X_{2})\stackrel{j}\to H^*(X_0)\to \cdots$$

From this sequence we get the short exact sequence $0\to \mathrm{im}\delta \to H^*(X)\to \mathrm{im}i\to 0$. By the exactness of this sequence, we have $\mathrm{im}i\cong \ker j$ and $\mathrm{im} \delta\cong H^{*-1}(X_0)/\ker \delta\cong H^{*-1}(X_0)/ \mathrm{im}j\cong \mathrm{coker} j$. This proves the lemma.

\begin{lemma}
Let $A$ be a generic $2\times 2$ Cartan matrix, then $H^*(BG(A))\cong \mathbb{Q}[\psi]$, where $\psi$ corresponds to Killing form which has degree 4. The Poincar\'{e} series of $BG(A)$ is $\D{\frac{1}{1-q^4}}$.
\end{lemma}

\noindent {\bf Proof: }By Lemma 2.6, we have the short exact sequence
$0\to \Sigma\mathrm{coker}j\to H^*(BG(A))\to \ker j\to 0$ with $j: H^*(BG_{\{1\}}(A))\oplus H^*(BG_{\{2\}}(A))\to H^*(BT)$.

The Poincar\'{e} series of $\textrm{coker}j$ is $\D{\frac{1}{(1-q^2)^2}-\frac{2}{(1-q^2)(1-q^4)}+\frac{1}{1-q^4}=0}$. Since $\ker j$ is isomorphic to $H^*(BT)^{W(A)}=\mathbb{Q}[\psi]$, its Poincar\'{e} series is $\D{\frac{1}{1-q^4}}$. Hence $H^*(BG(A))\cong \ker j\cong \mathbb{Q}[\psi]$.

This lemma is the special case of $n=2$ for Theorem 2.

\section{The proofs of the main theorems}
In this section we denote $BG_I(A)$ by $X_I$ for simplicity. For $I=\{i_1,i_2,\cdots,i_k\}$, we always denote $X_I$ by $X_{i_1i_2\cdots i_k}$. So we have $X_{\emptyset}=BT$ and $X_{12\cdots k} \simeq X_{12\cdots k-1}\cup_{ BT} X_k $.

\noindent{\bf Proof of Theorem 1:}
For $n=3$, $BG(A)$ is homotopic equivalent to $X_{12}\cup_{BT} X_{3}$. By Lemma 2.6 we have the short exact sequence
$0\to \Sigma \mathrm{coker}j\to H^*(BG(A))\to \ker j\to 0$. The homomorphism $j: H^*(X_{12})\oplus H^*(X_{3})\to H^*(BT)$ is given by $j(u,v)=Bi_{12}^*(u)-Bi_3^*(v)$, where $Bi^*_{12}$ and $Bi^*_3$ are induced by the homomorphisms $i_{12}: T\to G_{12}(A)$ and $i_3: T\to G_3(A)$. By Lemma 2.3, we observe that $\ker j$ is the sub-ring of Weyl group invariants. Since $A$ is non-symmetrizable, by Lemma 2.5, $\ker j\cong \mathbb{Q}$.
We have $\mathrm{im}j=\mathrm{im}(Bi^*_{12})+\mathrm{im}(Bi^*_3)$ and $\mathrm{im}(Bi^*_{12})\cap\mathrm{im}(Bi^*_3)=\ker j\cong\mathbb{Q}$. By Lemma 2.1, the Poincar\'{e} series of $\mathrm{im}(Bi^*_{12})$ and $\mathrm{im}(Bi_3)^*$ are $\D{\frac{1}{(1-q^2)(1-q^4)}}$ and $\D{\frac{1}{(1-q^2)^2(1-q^4)}}$ respectively. Combining these results, the Poincar\'{e} series of $\mathrm{coker} j$ is

$\D{\frac{1}{(1-q^2)^3}-\frac{1}{(1-q^2)(1-q^4)}-\frac{1}{(1-q^2)^2(1-q^4)}+1=\frac{2q^2-1}{(1-q^2)^2(1-q^4)}+1}$.

Hence for $n=3$ the Poincar\'{e} series of $H^*(BG(A))$ is $\D{q[\frac{2q^2-1}{(1-q^2)^2(1-q^4)}+1]+1}$.

For $n\geq 4$, we prove this theorem by induction on $n$. We assume the theorem is true for $n-1$. Since $A$ is non-symmetrizable, by Lemma 2.4, without loss of generality we can assume $A_{\{123\}}$ is non-symmetrizable.
Then $A'=A_{12\cdots n-1}$ is also non-symmetrizable. By Lemma 2.1, $H^*(X_{1,2\cdots n-1})\cong H^*(BG(A'))\otimes H^*(BS^1)$. By the induction assumption, the Poincar\'{e} series of $BG(A')$ is $\D{P_{n-1}(q)=q[\frac{(n-2)q^2-1}{(1-q^2)^{n-2}(1-q^4)}+1])+1}$. This means that the reduced cohomology $\widetilde H^*(BG(A'))$ concentrates in odd dimensions. Since $X\simeq X_{12\cdots n-1}\cup_{BT} X_n$, we use Lemma 2.6 to compute $H^*(BG(A))$. By the decomposition $H^*(X_{12\cdots n-1})\cong \widetilde H^*(BG(A'))\otimes H^*(BS^1)\oplus \mathbb{Q}\otimes H^*(BS^1)$, we have
$\ker j\cong \widetilde H^*(BG(A'))\otimes H^*(BS^1) \oplus \mathbb{Q}$. $\mathrm{im}j= \mathrm{im} Bi^*_{1,2\cdots n-1}+\mathrm{im} Bi^*_n$. And the intersection of $\mathrm{im} Bi^*_{1,2\cdots n-1}$ and $\mathrm{im} Bi^*_n $ is the sub-ring of Weyl group invariants. It is isomorphic to $\mathbb{Q}$.
The Poincar\'{e} series of $\mathrm{coker}j$ is $\D{\frac{1}{(1-q^2)^n}-\frac{1}{(1-q^2)^{n-1}(1-q^4)}-\frac{1}{1-q^4}+1}$. Therefore the Poincar\'{e} series of $BG(A)$ is
$\D{\frac{q}{(1-q^2)^n}-\frac{q}{(1-q^2)^{n-1}(1-q^4)}-\frac{q}{1-q^4}+q +\frac{q}{1-q^2}(P_{n-1}-1) +1}$
which is equal to $\D{q\left [\frac{(n-1)q^2-1}{(1-q^2)^{n-1} (1-q^4)} +1\right ]+1}$.

\noindent{\bf Proof of Theorem 2: }The proof of this theorem is similar to that of the Theorem 1. The difference is that for the symmetrizable case, the invariants of Weyl group is generated by the Killing form which is in degree 4.

We use induction on $n$. If $n=2$, by Lemma 2.7, the theorem is true. We assume this theorem is true for $n-1$. For an $n\times n$ symmetrizable generic Cartan matrix $A$, $A'=A_{12\cdots n-1}$ is a symmetrizable generic Cartan matrix. By induction assumption, the Poincar\'{e} series of $BG(A')$ is $\D{Q_{n-1}=\frac{1}{1-q^4}(q[\frac{(n-2)q^2-1}{(1-q^2)^{n-1}}+1]+1)}$. Since $BG(A)$ is homotopy equivalent to $X_{12\cdots n-1}\cup_{BT} X_n$. By a similar Mayer-Vietoris sequence computation, we get that the Poincar\'{e} series of $BG(A)$ is

\noindent $\D{Q_n=\frac{q}{(1-q^2)^n}-\frac{q}{(1-q^2)^{n-1}(1-q^4)}-\frac{q}{(1-q^2)(1-q^4)}+\frac{q}{(1-q^4)} +\frac{1}{(1-q^2)(1-q^4)}(Q_{n-1}-1)} $

\noindent $\D{+\frac{1}{1-q^4}=\frac{1}{1-q^4}\left [q(\frac{(n-1)q^2-1}{(1-q^2)^{n-1}} +1)+1\right ]}$.













\begin{remark}
In the proof of Theorem 2, we need the $H^*(BG(A))$-module structure on the Mayer-Vietoris sequences. In fact all the cohomology groups appear in the sequence are free $\mathbb{Q}[\psi]$-modules.
\end{remark}

\section{Some results derived from the main theorems}
In this section we need some general results in algebraic topology. For details see Whitehead\cite{Whitehead_78}.

From the expressions of the Poincar\'{e} series of $BG(A)$ in Theorem 1, we can see that for the non-symmetrizable case, the even dimensional cohomology group is $H^0(BG(A))$. Hence the cup product on $H^*(BG(A))$ is trivial and we have

\begin{coro}
For a generic $n\times n$ non-symmetrizable Cartan matrix $A$, the rational homotopy type of $BG(A)$ is $\bigvee\limits_{i=2}^\infty \bigvee\limits_{j=1}^{\alpha_{i}} S^{2i+1}$ with $P_n(q)=1 + \alpha_2 q^5+\alpha_3 q^7 +\cdots +\alpha_{i} q^{2i+1}+\cdots $.
\end{coro}

\begin{lemma}
For all $i\geq 2$, $\alpha_i=a_{2i}$.
\end{lemma}

\noindent{\bf Proof: }By definition we have
$$\D{\alpha_2 q^4+\alpha_3 q^6+\cdots=\frac{P_n(q)-1}{q}=\frac{(n-1)q^2-1}{(1-q^2)^{n-1} (1-q^4)} +1=1-\frac{1}{C_A(q)}=a_{4}q^4+a_6 q^6+\cdots+}$$
This proves the lemma.

By Milnor-Hilton theorem(see\cite{Whitehead_78}) and the homotopy equivalence $G(A)\simeq \Omega BG(A)$, we get

\begin{coro}
The homotopy Lie algebra $\pi_*(G(A))$ with Samelson product is the free graded Lie algebra generated by $\Sigma^{-1}\widetilde{H}_*(BG(A))$.

The Hopf algebra $H_*(G(A))$ is isomorphic to the tensor algebra $T(\Sigma^{-1} \widetilde{H}_*(BG(A)))$. 
\end{coro}

From the expressions of the Poincar\'{e} series in Theorem 2, we can see that for the symmetrizable case as $\mathbb{Q}[\psi]$-module the only even dimensional generator of $H^*(BG(A))$ is $1\in H^0(BG(A))$. The cup product on $H^*(BG(A))$ can be determined by the following lemma.

\begin{lemma}
Let $R=R_0\oplus R_1$ be a $Z_2$ graded ring and $R_1$ be a free $R_0$ module, if all the elements in $R_1$ are nilpotent, then $R_1 R_1=0$.
\end{lemma}

\noindent{\bf Proof: }Let $b_1,b_2\not=0$ be two elements in $R_1$ and $b_1b_2=a\in R_0$. We show $a=0$. Since $b_2$ is nilpotent, there exists integer $k>0$ such that $b_2^{k-1}\not=0$ but $b_2^k=0$. Then we have $a b_2^{k-1} =b_1 b_2^k=0$. Hence we get $a=0$ from the fact that $R_1$ is a free $R_0$ module.

Set $R_0=H^{even}(BG(A))$ and $R_1=H^{odd}(BG(A))$, by this lemma, we get

\begin{coro}
For a generic $n\times n$ symmetrizable Cartan matrix $A$, the rational homotopy type of $BG(A)$ is  $BS^3\times\bigvee\limits_{i=2}^\infty \bigvee\limits_{j=1}^{\beta_{i}} S^{2i+1}$ with $(1-q^4)Q_n=1 + \beta_2 q^5+\beta_3 q^7 +\cdots +\beta_{i} q^{2i+1}+\cdots $.
\end{coro}

Similarly we have

\begin{lemma}
For all $i\geq 2$, $\beta_i=b_{2i}$.
\end{lemma}

\begin{coro}
The homotopy Lie algebra $\pi_*(G(A))$ with Samelson product is the direct sum of $\pi_*(S^3)$ and the free graded Lie algebra generated by $\Sigma^{-1}\bar{H}_*(BG(A))$, where $\bar H_*(BG(A))\cong \widetilde H_*(\bigvee\limits_{i=2}^\infty \bigvee\limits_{j=1}^{\beta_{i}} S^{2i+1}) $.

The Hopf algebra $H_*(G(A))$ is isomorphic to the $\mathbb{Q}$-algebra $H_*(S^3)\times T(\Sigma^{-1} \bar{H}_*(BG(A)))$. 
\end{coro}

Now we have the following result.

\begin{prop}
If $A_1,A_2$ are two generic Cartan matrices of size $n_1$ and $n_2$, then $G(A_1)$ and $G(A_2)$ are rational homotopy equivalent Hopf spaces if and only if $n_1=n_2$ and $\epsilon(A_1)=\epsilon(A_2)$.
\end{prop}

\end{document}